\font\tenmsb=msbm10
\font\sevenmsb=msbm7
\font\fivemsb=msbm5
\newfam\msbfam
  \textfont\msbfam=\tenmsb
  \scriptfont\msbfam=\sevenmsb
  \scriptscriptfont\msbfam=\fivemsb
\def\Bbb{\fam\msbfam\tenmsb}

\newcommand{\ext}[1]{\stackrel{#1}{\wedge}}
\newcommand{\Caut}{\mbox{${{\rm Caut}(\pi_1(X))}$}}
\newcommand{\Vaut}{\mbox{${{\rm Vaut}(\pi_1(X))}$}}
\newcommand{\Vautt}{\mbox{${{\rm Vaut}(\pi_1({\tilde X}))}$}}

\def\VautGa{{\rm Vaut}(\Gamma)}
\def\Hqs{{\rm Homeo_{q.s.}}(S^1)} 
\def\Mob{{\rm Mobius}(S^1)} 

\def\ra{\rightarrow }

\def\Q{{\Bbb Q}}

\def\C{{\Bbb C}}
\def\C*{{\C}^{*}} 
 
\def\Z{{\Bbb Z}}

\def\Cg{{\cal C}_g}

\def\Hin{{H}_{\infty}}
\def\HinX{{H}_{\infty}(X)}

\def\Tg{{\cal T}_g}
\def\Th{{\cal T}_h}

\def\Tin{{\cal T}_{\infty}}
\def\TinX{{\cal T}_{\infty}(X)}
\def\TchinX{{\cal T}_{\infty}^{ch}(X)}
\def\TinY{{\cal T}_{\infty}(Y)}
\def\THin{{\cal T}(\Hin)} 
\def\THinX{{\cal T}(\Hin(X))}

\def\Mg{{\cal M}_g}

\def\Mt{{\cal M}_{\tilde g}}

\def\Min{{\cal M}_{\infty}(X)}
\def\MinX{{\cal M}_{\infty}(X)}

\def\MCin{{MC}_{\infty}}

\def\Kch{{K}^{ch}(X)}

\def\D0Q{{DET}(0,{\Q})}

\def\Tpi{{\cal T}(\pi)}
\def\Mpi{{\cal M}(\pi)}

\def\tg{{\tilde g}}
\def\tX{{\tilde X}}

\def\tom{{\tilde \omega }}
\def\lt{{\tilde{\lambda}}}

\def\A{{\cal A}}

\def\C{{\cal C}}

\def\M{{\cal M}}
\def\O{{\cal O}}

\def\T{{\cal T}}

\def\a{\alpha }
\def\b{\beta }
\def\D{\Delta }
\def\ga{\gamma }
\def\Ga{\Gamma }

\def\r{\rho }
\def\p{\psi }

\def\l{\lambda }
\def\La{\Lambda }

\def\s{\sigma }
\def\Si{\Sigma }

\def\om{\omega }

\def\VautGa{{\rm Vaut}(\Gamma)}
\def\Hqs{{\rm Homeo_{q.s.}}(S^1)} 

\def\Pa{{\rm Pic_{alg}}(\Mg)}

\def\Pf{{\rm Pic_{fun}}(\Mg)}

\long\def\comment#1\endcomment{}

\def\mapright#1{\smash{
    \mathop{\longrightarrow}\limits^{#1}}}

\def\mapdown#1{\Big\downarrow   
 \rlap{$\vcenter{\hbox{$\scriptstyle#1$}}$}}
\documentstyle[12pt]{article}
\begin{document}
\baselineskip=16pt
\begin{flushright}
{\it To appear in:} {\underline {{\bf Contemporary Math. series,} 
({\it Amer. Math. Soc.,})}}\\ 
{\underline {{\it Bers Colloquium (New York, 1995) volume.}}}\\ 
\end{flushright}

\begin{center}
{\bf WEIL-PETERSSON GEOMETRY AND DETERMINANT BUNDLES}\\ 
{\bf ON INDUCTIVE LIMITS OF MODULI SPACES}
\footnote{
~~~~~~~~~~~~Preprint no.: imsc-96/05/15}

\vspace{.4cm}

{\bf  Indranil Biswas ~ {\it and} ~ Subhashis Nag}\\

\vspace{.6cm}

{\small {\bf Abstract}} 
\end{center}

\baselineskip=12pt
In the paper [BNS] the authors and Dennis Sullivan 
constructed the universal direct system of
the classical Teichm\"uller spaces of Riemann surfaces of varying genus.
The direct limit, which we called the universal commensurability 
Teichm\"uller space, $\Tin$, was shown to carry on it a natural 
action of the universal commensurability mapping class group, $\MCin$. 
In this paper we identify an interesting cofinal sub-system corresponding 
to the tower of finite-sheeted {\it characteristic} coverings over 
any fixed base surface $X$. Utilizing a certain subgroup $\Caut$ 
inside $\MCin$, (associated intimately to this characteristic tower), we 
descend to an inductive system of {\it moduli} spaces, and construct 
the direct limit ind-variety $\Min$.

Invoking curvature properties of Quillen metrics on determinant bundles,
and naturality under finite coverings of Weil-Petersson forms, we are
able to construct on $\Min$ the natural sequence of determinant of
cohomology line bundles, as well as the Mumford isomorphisms connecting
these. 

\baselineskip=16pt
\vspace{1cm}
\noindent 
{\small {\bf I. INTRODUCTION}} 
\medskip

Let $\Tg$ denote the Teichm\"uller space comprising compact marked
Riemann surfaces of genus $g$, and $\Mg$ be the moduli
space of Riemann surfaces of genus $g$ obtained by quotienting $\Tg$ by
the action of the mapping class (=modular)  group, $MC_g$. Denote
by $DET_n \rightarrow \Tg$ the line bundle given by the
determinant of cohomology construction for the $n$-th
tensor power ($n \in \Z$) of the relative cotangent bundle on 
the universal family of Riemann surfaces, $\Cg$, over $\Tg$;
(see IV.1 below and [BNS], [D], for detailed definitions).
The bundle $DET_0$ is classically called the Hodge line bundle;
it is a fundamental fact that Hodge generates the entire Picard
group of the moduli functor (see the section cited above).

Each bundle $DET_n$ comes equipped with a hermitian structure 
which is obtained from the construction of Quillen of metrics on
determinant bundles, [Q]. Quillen's construction is subordinate to the 
choice of a smoothly varying family of K\"ahler metrics on the fibers of 
the family of Riemann surfaces; we utilize the Poincar\'e hyperbolic metric 
on the fibers of $\Cg$ (for $g\geq 2$) to obtain the corresponding natural 
Quillen metric on each $DET_n$.

By applying the Grothendieck-Riemann-Roch theorem, Mumford [Mum] had 
shown that $DET_n$ is a certain fixed (genus-independent) tensor power 
of the Hodge bundle over each moduli space $\Mg$. Precisely:
$$
DET_{n} ~=~ DET_{0}^{\otimes (6n^2-6n+1)}
\leqno{(1.1)}
$$
The isomorphism may be considered as an equivariant isomorphism of
$MC_g$ equivariant line bundles over $\Tg$. The Mumford
isomorphism is unique up to a non-zero multiplicative constant, and
can be chosen to be an isometry with respect to the Quillen metrics
mentioned above.

There is a very interesting connection, discovered by Belavin and Knizhnik
[BK], between the Mumford isomorphism above for the case $n=2$, (namely 
that $DET_{2}$ is the $13$-th tensor power of the Hodge bundle), and the 
existence of the Polyakov string measure on the moduli space $\Mg$. 
For an exposition of this connection see, for instance, [N2].
That suggests the natural question of finding a {\it genus-independent} 
formulation of the Mumford isomorphisms over some ``universal'' 
parameter space of Riemann surfaces (of varying genus). 

Our joint paper with Dennis Sullivan gives such a genus-independent,  
universal version of the determinant bundles and Mumford's isomorphism
by working over the universal commensurability Teichm\"uller space.
The geometrical objects in [BNS] exist over this universal base space 
$\Tin = \TinX$, which is defined as the infinite direct limit of the 
Teichm\"uller spaces of higher genus pointed surfaces that are 
finite unbranched coverings of any pointed reference surface $X$. The 
bundles and the relating isomorphisms are equivariant with respect 
to the natural action of a large new mapping class group, called
the  universal commensurability group $\MCin$ -- which we introduced in
[BNS]. Our method there was to utilize a subtle form of the
Grothendieck-Riemann-Roch theorem in a formulation of Deligne, [D], which
depended on a certain construction of Deligne known as the ``Deligne
pairing''.

The main purpose of the present paper is to obtain a genus-independent
description of the Mumford isomorphisms over {\it inductive limits of 
moduli spaces} $\Mg$ by looking at the inverse system of finite 
unbranched {\it characteristic} coverings of any reference surface $X$. The
characteristic covers are shown to form a cofinal tower (in the tower of
all finite unbranched coverings of $X$), and the construction proceeds
over the direct limit of {\it moduli} spaces (rather than at the
Teichm\"uller level). We consequently obtain certain ``rational line 
bundles'' over the direct limit, $\Min$, of moduli spaces, with their 
relating Mumford isomorphisms.  We investigate the relationship between
$\Min$ and $\Tin$ by considering the subgroup of the universal
commensurability modular group that acts on $\Tin$ to produce $\Min$ as
the quotient. The representation of the commensurability modular group as
a subgroup of the group of quasisymmetric homeomorphisms of the circle,
and the relation with the classical Teichm\"uller theory of the
Ahlfors-Bers universal Teichm\"uller space, are also explained here. 
We present the material in a more leisurely fashion than in [BNS], 
highlighting also some salient questions that remain unresolved.

Another purpose of this article is to show that one can use the 
Weil-Petersson K\"ahler geometry of the Teichm\"uller spaces 
to obtain the desired genus-independent construction of $DET$ bundles 
and Mumford isomorphisms in some special but interesting cases,
instead of the more sophisticated GRR theorem invoked in [BNS], where we
worked in a very general set-up. The Weil-Petersson form comes into 
play because it represents (up to scaling factors) the curvature form 
for each of the $DET_n$ bundles (when these bundles are equipped with 
their Poincar\'e-Quillen metrics). The formulation of our final results 
turns out to be somewhat different from the theorem we presented in [BNS].  

The parameter spaces obtained by passing to the direct limit of the
Teichm\"uller or moduli spaces over varying genus, can be interpreted as 
a certain space of (``transversely locally constant'') complex structures
on the corresponding {\it solenoidal surface} arising by taking the
inverse limit (through the tower of coverings) of the classical compact 
surfaces. There is an interplay between the topological type of the 
solenoidal inverse limit and the type of the associated direct limit 
moduli space, which also appears in the work presented in this paper.
 
\medskip
\noindent
{\it Acknowledgement:} We would like to thank Dennis Sullivan for 
many discussions, and for suggesting to us the idea of utilizing 
characteristic coverings.

\bigskip
\noindent
{\small {\bf II. THE UNIVERSAL DIRECT LIMIT $\Tin$}}
\medskip

\noindent
{\bf II.1. Coverings and the Teichm\"uller functor:}
We start with a fundamental topological situation. Let
$$
\pi~ : ~ \tX ~ \longrightarrow ~ X
\leqno(2.1)
$$ 
be an {\it unramified finite covering}, orientation preserving, 
between two compact connected oriented two manifolds $\tX$ and $X$ of 
genera $\tg$ and $g$, respectively. Assume $g \geq 2$. The degree of the
covering $\pi$, which will play an important role, is the ratio
of the respective Euler characteristics; namely, $\deg (\pi)=(\tg-1)/(g-1)$.

The Teichm\"uller space ${\cal T}_g$ (resp. ${\cal T}_{\tg}$) is the quotient
of all complex structure on $X$ (resp. $\tX$) by the 
group of all orientation preserving
diffeomorphisms of $X$ (resp. $\tX$) which are homotopic to the identity
map -- this group will be denoted by ${\rm Diff}^+_0(X)$ (resp.
${\rm Diff}^+_0(\tX)$).
Given any complex structure on $X$, we may pull back this structure
via $\pi$ to get a complex structure on $\tX$. The homotopy lifting 
property guarantees that there is a diffeomorphism 
${\tilde f}\in{\rm Diff}_0({\tX})$ 
which is a lift of any given $f \in{\rm Diff}_0({X})$. 
It follows that the process of pulling back complex structure from 
$X$ onto $\tX$ induces a well-defined map at the level of the 
Teichm\"uller spaces:
$$
{\cal T}(\pi): ~ {\cal T}_g ~ \longrightarrow ~ {\cal T}_{\tg}
\leqno(2.2)
$$ 
It is known that this map ${\cal T}(\pi )$ is a {\it proper holomorphic 
embedding} between these finite dimensional complex manifolds;
furthermore, the map $\Tpi$ respects the quasiconformal-distortion 
(=Teichm\"uller) metrics. 

Two coverings are said to be in the same homotopy class if they
are homotopic through continuous mappings. When working with 
pointed surfaces and base-point preserving coverings, we shall
say that two coverings are in the same based homotopy class if they
are homotopic through a base-point preserving family of continuous 
mappings. It is easy to see that the above embedding between the
Teichm\"uller spaces depends only on the unbased homotopy class of
the covering $\pi$. (The pullback of a given complex structure on $X$ to
$\tX$, using a covering $\pi$, depends of course on the map itself, 
and not just on its homotopy class. But this dependence disappears 
when passing to the level of the corresponding Teichm\"uller spaces.)

At the level of Fuchsian groups, one should note that the covering
space $\pi$ corresponds to the choice of a subgroup $H$ of finite index 
(=$\deg(\pi)$) in the uniformizing group $G$ for $X$, and the embedding 
(2.2) is then the standard inclusion mapping $\T(G) \rightarrow \T(H)$; 
(see Chapter 2, [N1]).

\medskip
\noindent
{\bf Remark 2.3:}~ One notices that the morphisms of the type ${\cal 
T}(\pi)$ in (2.2) constitute a contravariant functor from the category
whose objects are closed oriented topological surfaces 
and the morphisms being the covering maps, to the 
category of finite dimensional complex manifolds and the holomorphic 
embeddings. This functor will be denoted by $\cal T$. We 
shall have more to say along these lines below. 
\medskip

We construct a category $\A$ of certain topological objects and morphisms: 
the objects, $Ob(\A)$, constitute a set of compact oriented topological 
surfaces each equipped with a base point ($\star$), there being exactly 
one surface of each genus $g \geq 0$; let the object of genus $g$ be
denoted by $X_g$. The morphisms are based homotopy classes of pointed 
covering mappings
$$
\pi ~ : ~ (X_\tg, \star) ~ \longrightarrow ~ (X_g, \star)
$$ 
there being one arrow for each such based homotopy class. 
An important point to note is that the
monomorphism of fundamental groups induced by any representative of
the based homotopy class of coverings $\pi$ is unambiguously defined. 

\medskip
\noindent
{\bf II.2. The direct system of classical Teichm\"uller spaces:}
Fix a genus $g$ and let $X = X_g$. Observe that all the morphisms 
with the fixed target $X_g$: 
$$
K_g ~ = ~ K(X) ~ = ~ \{\a \in {\rm Mor}(\A) : {\rm Range}(\a)=X \}
\leqno(2.4)
$$
constitute a {\it directed set} under the partial ordering given by 
factorization of covering maps. Thus if $\a$ and $\b$ are two morphisms 
in the above set, then $\b \succ \a$ if and only if the image of the
monomorphism $\pi_1(\b)$ is contained within the image of $\pi_1(\a)$.
This happens if and only if there is a commuting triangle of morphisms:
$\b = \a \circ \theta$. It is important to note that the factoring 
morphism $\theta$ is {\it uniquely} determined because we are working 
with surfaces with base points. 

\medskip
\noindent
{\it Remark:} Notice that the object of genus $1$ in $\A$ only has 
morphisms to itself -- so that this object together with all its 
morphisms (to and from) form a subcategory.
\medskip

As shown in (2.2), each morphism of $\A$ induces a proper, holomorphic, 
Teichm\"uller-metric preserving embedding between the corresponding 
finite-dimensional Teichm\"uller spaces. We can thus create the natural 
{\it direct system of Teichm\"uller spaces} over the above directed 
set $K_g$, by associating to each $\a \in K_g$ the Teichm\"uller space 
$\T(X_{g(\a)})$, where $X_{g(\a)} \in Ob(\A)$ denotes the domain surface
for the covering $\a$. To each $\b \succ \a$ one associates the 
corresponding holomorphic embedding $\T(\theta)$ (with
$\theta$ as above). From this direct system 
we form the {\it direct limit Teichm\"uller space over $X=X_g$}:
$$
\Tin(X_g) = \TinX := {\rm ind. lim.} \T(X_{g(\a)})
\leqno(2.5)
$$
This limit $\TinX$ is an ``ind-space'' in the sense of 
Shafarevich [Sha]. In other words, it an inductive limit of 
finite dimensional spaces. It is a metric space with a well-defined 
Teichm\"uller metric. Indeed, $\TinX$ also carries a natural 
Weil-Petersson Riemannian structure obtained from scaling the 
Weil-Petersson pairing on each finite dimensional stratum, $\Th$, 
by the factor $(h-1)^{-1}$. In fact, compare Theorem 9.1 of [NS]
(asserting the existence of Weil-Petersson structure on the Teichm\"uller 
space $\THinX$) with the crucial Lemma 5.1 below.

The space $\Tin$ is called the {\it universal commensurability Teichm\"uller 
space}: it is an universal parameter space for compact Riemann surfaces. 
$\Tin$ serves as the base space for our construction of universal Mumford 
isomorphisms in [BNS].

\medskip
\noindent
{\bf II.3. The Teichm\"uller space, $\THin$, of the hyperbolic solenoid:} 
Over the very same directed set $K_g$ in (2.4), we may also define a
natural {\it inverse system of surfaces}. This is done by associating to
each $\a \in K_g$ a certain copy, $S_{\a}$ of the pointed surface 
$X_{g(\a)}$.
[Note: Fix a universal covering over of $X=X_g$. The surface $S_{\a}$
can be taken to be this universal covering quotiented 
by the action of the subgroup $Im(\pi_1(\a)) \subset {\pi_1}(X,\star)$ 
using the action of the deck transformations.] If $g \ge 2$, 
then the inverse limit of this system is the {\it universal solenoidal 
surface} $H_{\infty}(X) = {\rm inv~lim}X_{g(\a)}$, that was studied in 
[S],[NS]. 

The ``universality'' of this object resides in the evident but important 
fact that these spaces $\HinX$, as well as their 
Teichm\"uller spaces $\THinX$, 
do {\it not} really depend on the choice of the base surface $X$. If we 
were to start with a surface $X'$ of different genus (both genera 
being greater than one), we could pass to a common covering
surface of $X$ and $X'$ (always available!), and hence the limit 
spaces we construct would be naturally isomorphic. We are therefore 
justified in suppressing $X$ in our notation and referring to $\HinX$ 
as simply $\Hin$.

The space $H_{\infty}$ is compact. For each surface $X$
(of genus greater than one) there is a natural fibration
$$
\pi_{\infty}(X): ~ H_{\infty} \rightarrow  X,
$$
the fibers being Cantor sets. The path components of $H_{\infty}$
are called ``leaves". Each leaf, with the ``leaf-topology" it
inherits from $H_{\infty}$, is a simply connected two-manifold,
and the restriction of $\pi_{\infty}(X)$ to any leaf is a
universal covering of $X$. There are uncountably many leaves in
$H_{\infty}$, and each is a dense subset of $H_{\infty}$. 
Each leaf is thus identifiable with a hyperbolic plane. 
That is why we call $\Hin$ the universal hyperbolic solenoid. 
The facts above follow from a careful study of this inverse system of
surfaces, the main tool being the lifting of paths in $X$ to its coverings.

As explained in [S],[NS], the solenoid $\Hin$ has a natural
Teichm\"uller space comprising equivalence classes of complex structures 
on the leaves -- the leaf complex structures being required to vary 
continuously in the fiber (Cantor) directions. In particular, any
complex structure assigned to any of the surfaces $X_{g(\a)}$
appearing in the inverse tower can be pulled back to all the surfaces
above it -- and therefore assigns a complex structure of the sort
demanded on $\Hin$ itself. These complex structures that arise from
some finite stage can be characterized as the ``transversely locally
constant'' (TLC) ones (see [NS]), and they comprise precisely the 
{\it dense} subset $\TinX$ sitting within the separable Banach manifold 
$\THinX$. We collect the above discussions in the:

\medskip
\noindent
{\bf Proposition 2.6 [BNS]:} {\it The ind-space $\TinX$ arises as  
an inductive limit of finite dimensional complex manifolds, and hence
carries a complex structure defined strata-wise. The completion of 
$\TinX$ with respect to the Teichm\"uller metric is the separable 
complex Banach manifold $\THinX$. 

In fact, $\TinX$ can be embedded in Bers' universal Teichm\"uller 
space, $\T(\Delta)$, ($\Delta$ denotes the unit disc), as a directed 
union of the Teichm\"uller spaces of a family of Fuchsian groups. The 
Fuchsian groups vary over the finite index subgroups of a fixed 
cocompact Fuchsian group $G$, $X=\Delta/G$. The closure in $\T(\Delta)$ 
of $\TinX$ is a Bers-embedded copy of $\THinX$.}

\medskip
\noindent
{\bf II.4. The commensurability mapping class group $\MCin$:}
We proceed to recall in some detail a construction introduced in
[BNS]. A remarkable fact about the situation above is that every 
morphism $\pi:Y \longrightarrow X$ of $\A$ induces a natural Teichm\"uller 
metric preserving {\it homeomorphism} 
$$
\Tin(\pi) ~ : ~ \TinY ~ \longrightarrow ~ \TinX
$$
The map $\Tin(\pi)$ is invertible simply because the morphisms of $\A$ 
with target $Y$ are cofinal with those having target $X$ (thus all 
finite ambiguities are forgotten in passing to the inductive limits!). 
It is also clear that $\Tin(\pi)$ is a biholomorphic identification 
(with respect to the strata-wise complex structures). Recall the
functor $\cal T$ defined in Remark 2.3. We may similarly define a
functor using the morphisms $\Tin (\pi)$, which will be denoted
by $\Tin$. Note that the functor $\Tin$ is covariant -- whereas
the Teichm\"uller functor $\T$ itself was contravariant.

For a given pair of coverings (not necessarily homotopic)
$$
\a: Y \longrightarrow X
\hspace{.5in} {\rm and} \hspace{.5in}
\b: Y \longrightarrow X
\leqno{(2.7)}
$$
we have an automorphism
$${\cal T}_{\infty}({\b})\circ {\cal T}_{\infty}(\a)^{-1}
\leqno{(2.8)}
$$
of ${\cal T}_{\infty}(X)$. This automorphism preserves the metric
on $\TinX$ and hence it extends to the metric completion of it.

We will call a pair of the form (2.7) a finite {\it self correspondence}
of $X$.

More generally, assume that we are given a cycle of coverings starting
and ending at $X$:
$$\matrix
{Y_{k} & \hbox{---} & Y_{k+1}\cr
 \vert && \vert \cr
Y_{k-1} &  & Y_{k+2}\cr
\vert && \vert \cr
\vdots &  & \vdots \cr
Y_1 &  & Y_{n} \cr
\vert &&  \vert \cr
X & = & X\cr} \leqno{(2.9)}$$
where $X$, $Y_i$ are all objects of the category $\cal A$ and
all horizontal and vertical lines represent morphisms (pointing in
arbitrary directions) of $\cal A$. Using the automorphism in (2.5) 
for each covering in the diagram, and applying it to all
the coverings in (2.9), we get an automorphism of $\TinX$ just as in (2.8).
Note that since ${\cal T}_{\infty}(\pi)$ in (2.5) is invertible,
the horizontal and the vertical lines in (2.9) are allowed to be maps
in any direction. For example, if some of the maps
$Y_i \, \hbox{---} \, Y_{i-1}$ 
point upwards and some downwards, or left/right, in any such instance 
the construction of the automorphism of $\TinX$ 
(obtained by following the entire cycle around) remains valid.

Thus we see that each $\TinX$, and consequently also 
its metric completion $\THinX$, is equipped with a large {\it 
automorphism group} -- one from each such undirected cycle of 
morphisms of $\A$ starting from $X$ and returning to $X$. 
By repeatedly using pull-back diagrams (i.e., by choosing the appropriate 
connected component of the fiber product of covering maps), it is fairly
easy to see that the automorphism of $\TinX$ arising from
any (many arrows) cycle can be obtained simply from a
self-correspondence, i.e., a two-arrow cycle.

These self-maps constitute a group of biholomorphic automorphisms 
of $\TinX$ that we shall call the {\it universal commensurability modular 
group} $\MCin(X)$, acting on $\TinX$ and on $\THinX$. 
We shall show below (Proposition 2.17) how $\MCin(X)$ may be 
realized as a subgroup of the classical universal modular group. 

To clarify matters further, we consider the abstract graph 
($1$-complex), $\Ga(\A)$, obtained from the topological 
category $\A$ by looking at the objects as vertices and 
the (undirected) arrows as edges. It is clear from the definition 
above that the fundamental group of this graph, viz. $\pi_{1}(\Ga(\A),X)$, 
is acting on $\TinX$ as these automorphisms. We may fill in all triangular 
$2$-cells in this abstract graph whenever two morphisms (edges) compose to 
give a third edge; the thereby-reduced fundamental group of this
$2$-complex can be shown to produce faithfully the action of $\MCin(X)$
on $\TinX$. 

\noindent 
{\it Remark on the genus one subcategory:}  
For the genus one object $X_1$ in $\A$, we can make the entire business
explicit. We know that the Teichm\"uller space for any unramified 
covering is a copy of the upper half-plane $H$. The maps $\Tpi$ are
M\"obius identifications of copies of the half-plane with itself, 
and we easily see that the pair $(\Tin(X_1),\MCin(X_1))$ is 
identifiable as $(H,PGL(2,\Q))$. Notice that the action has dense
orbits in this case. Anticipating for a moment the definition of
the virtual automorphism group, ${\rm Vaut}$, given in II.5 below, we
remark that $GL(2,\Q)$ is indeed 
${\rm Vaut}(\Z \oplus \Z)$, and ${\rm Vaut}^{+}$ is the subgroup of index 
$2$ therein, as expected.

In the general case, if $X \in Ob(\A)$ is of any genus $g \geq 2$,
then we get an infinite dimensional ``ind-space'' as $\TinX$ with  
the action of $\MCin(X)$ on it as described. Since the tower of coverings
over $X$ and $Y$ (both of genus higher than $1$) eventually become 
cofinal, it is clear that {\it for any choice of genus higher than one 
we get {\bf one isomorphism class} of pairs} $(\Tin, \MCin)$. 

\noindent
{\bf II.5. Virtual automorphism group of $\pi_{1}(X)$ and $\MCin$:}
In the classical situation, the action of the mapping class group 
$MC(X)$ on $\T(X)$ was induced by the action of (homotopy classes of) 
self-homeomorphisms of $X$; in the direct limit set up we  now have the
more general (homotopy classes of) self-correspondences of $X$ inducing 
the new mapping class automorphisms on $\TinX$. In fact, we will see that
our group $\MCin$ corresponds to ``virtual automorphisms'' of the 
fundamental group $\pi_{1}(X)$, -- generalizing exactly the classical
situation where the usual $Aut(\pi_{1}(X))$ appears as the action 
via modular automorphisms on $\T(X)$.

Given any group $G$, one may look at its ``partial'' or ``virtual'' 
automorphisms, [Ma]; as opposed to usual automorphisms which are defined 
on all of $G$, for virtual automorphisms we demand only that they 
be defined on some finite index subgroup. To be precise, consider
all isomorphisms $\r :H \longrightarrow K$ where $H$ and $K$ are
subgroups of finite index in $G$. Two such isomorphisms (say $\r_1$ 
and $\r_2$) are considered equivalent if there is a finite index subgroup
(sitting in the intersection of the two domain groups) on which they
coincide. The equivalence class $[\r]$  -- which is like the {\it
germ} of the isomorphism $\r$ -- is called a {\it virtual automorphism} 
of $G$; clearly the virtual automorphisms of $G$ constitute a group, 
christened ${\rm Vaut}(G)$, under the obvious law of composition, (i.e., 
compose after passing to deeper finite index subgroups, if necessary).  

Clearly ${\rm Vaut}(G)$ is trivial unless $G$ is infinite (though there do
exist infinite groups -- see [MT] -- such that ${\rm Vaut}$ is trivial).  
Also evident is the fact that $${\rm Vaut}(G) ~ = ~ {\rm Vaut}(H)$$
where $H$ is a finite index subgroup of $G$. Since we shall apply this 
concept of virtual automorphism to the fundamental group of a surface of
genus $g$, ($g>1$), the last remark shows that our 
${\rm Vaut}(\pi_{1}(X_g))$ {\it is genus independent}! 

In fact, ${\rm Vaut}$ presents us a neat way of
formalizing the ``two-arrow 
cycles'' (2.7) which we introduced to represent elements of $\MCin$. 
Letting $G = \pi_{1}(X)$, (recall that $X$ is already equipped with a 
base point), we see that the diagram (2.7) corresponds exactly to the 
following virtual automorphism of $G$:
$$
[\r] ~ = ~ [{\b}_{*}\circ{\a}_{*}^{-1}:{\a}_{*}(\pi_{1}(\tX)) 
\longrightarrow {\b}_{*}(\pi_{1}(\tX))]
$$
Here ${\a}_{*}$ denotes the monomorphism of the fundamental group 
$\pi_{1}(\tX)$ into $\pi_{1}(X) = G$ induced by $\a$, 
and similarly ${\b}_{*}$ etc..
We let ${\rm Vaut}^{+}({\pi}_{1}(X))$ denote the subgroup of
${\rm Vaut}$ arising from pairs of {\it orientation preserving} coverings.
(We shall ignore the difference
between $\Vaut$ and ${\rm Vaut}^{+}({\pi}_1(X))$ below -- when speaking
of ${\rm Vaut}$ we shall mean the ${\rm Vaut}^{+}$.) 

\noindent
{\it Remark:}
The reduction of any many-arrow cycle in $\Ga(\A)$ to a two-arrow cycle
utilizes successive fiber product diagrams; there is some amount 
of choice in this reduction process, and one may obtain different 
two-arrow cycles starting from the same cycle; however, one may verify
that the virtual automorphism that is defined via any reduction is 
unambiguous. 

The final upshot is:

\noindent
{\bf Proposition 2.10:} 
{\it (a) ${\rm Vaut}^{+}({\pi}_{1}(X))$ is naturally isomorphic to 
$\MCin(X)$. 

(b) The natural homomorphism:
$ {\pi}_{1}(\Ga(\A)_{fill},X) \rightarrow {\rm Vaut}^{+}({\pi}_{1}(X))$ 
is an isomorphism. Here $\Ga(\A)_{fill}$ denotes the $2$-complex
obtained from the graph $\Ga(\A)$ by filling in all commuting triangles in
$\Ga(\A)$.}

\noindent
{\it Summarizing remark:}
So, interestingly enough, the usual $Aut({\pi}_{1}(\tX))$ acts as the
standard modular action on each of the classical Teichm\"uller spaces,
$\T(\tX)$, which constitute the various finite dimensional 
strata in $\Tin(X)$ (associated to surfaces of varying genus),
 -- whereas the direct limit Teichm\"uller space is acted upon by 
this (genus-independent) new modular group $\Vaut = \Vautt$. 

\medskip
\noindent
{\bf II.6. Representation of $\Vaut$ within $Homeo(S^1)$:}
${\rm Vaut}(\pi_{1}(X))$ allows certain natural representations in the 
homeomorphism group of the unit circle $S^1$, by the standard theory of
{\it boundary homeomorphisms} (see, for example, Chapter 2, [N1]). In fact, 
we get one such representation for each choice of cocompact Fuchsian 
group $\Ga$ faithfully representing $\pi_1(X)$. We take the base point on 
$X$ to be the image of the origin of the unit disc under the universal 
covering projection $u:\D \ra \D/{\Ga} \equiv X$.

Thus let $[\r] \in {\rm Vaut}(\Ga)$ be represented by the isomorphism 
$\r:H \ra K$. Then the Fuchsian subgroups $H$ and $K$ represent,
respectively, the (pointed) Riemann surfaces $Y=\D/H$ and $Z=\D/K$ 
covering $X$. The base points on $Y$ and $Z$ are, of course, 
the respective images of the origin of $\D$. The given isomorphism 
$\r: \pi_1(Y) \ra \pi_1(Z)$ can now be realized (using Nielsen's theorem) 
by an orientation preserving diffeomorphism
(quasiconformal homeomorphism is enough for our purposes) 
$h_\r:Y \ra Z$, preserving base points, satisfying: 
$$
\pi_{1}(h_{\r})~ = ~ \r~ : ~ H ~ \longrightarrow ~ K
$$   
The based homotopy class of $h_\r$ is uniquely determined. We {\it lift}
$h_\r$ to the universal covering to get a self-diffeomorphism ${\Si}_\r$ 
of $\D$ preserving the origin.

$$
\matrix{
{\D}
&\mapright{\Si_\r}
&{\D}
\cr 
\mapdown{}
&
&\mapdown{}
\cr 
Y
&\mapright{h_\r}
&Z 
\cr}
$$ 

The basic equation relating $\Si_\r$ to $\r$ is:
$$
\Si_{\r} \circ h \circ {\Si_\r}^{-1} = \r(h), ~~for~all~h \in H.
\leqno(2.11)
$$

Now associate to $[\r] \in {\rm Vaut}(\Ga)$ {\it the boundary values
of this lift of} $h_\r$ to obtain the desired representation:
$$
\Si ~ : ~ {\rm Vaut}(\Ga) ~ \longrightarrow ~{\rm Homeo_{q.s.}}(S^1);~~ 
\Si([\r]) ~ = ~ \partial \Si_{\r} ~ := ~ {\Si}_{\r}\vert_{\partial\Delta} 
\leqno(2.12) 
$$
Since we are dealing with compact surfaces, any diffeomorphism is
quasiconformal -- hence so is the lift $\Si_{\r}$. The boundary homeomorphism 
$\partial \Si_{\r}$ therefore exists by continuous extension, and is 
a quasisymmetric homeomorphism on $S^1 = \partial\D$. That boundary 
homeomorphism depends only on the homotopy class of $h_\r$ for well-known
reasons -- see, for example, pp. 114ff of Chapter 2 of [N1].

Consequently, (2.12) can be seen to be well-defined on equivalence classes
$[\r]$, and it is not hard to check that indeed $\Si$ gives us a 
{\it faithful} representation  of $\VautGa$ within $\Hqs$.

\noindent
{\it A simple description of the boundary homeomorphism:} Given the
virtual automorphism $\r:H \ra K$, consider the natural map it defines of the
orbit of the origin (=$0$) under $H$ to the orbit of $0$ under $K$. Namely:
$$
\s_\r ~ : H(0) \longrightarrow K(0); 
\hspace{.5in} 
h(0) ~ \longmapsto ~ \r(h)(0) 
$$
But each orbit under these cocompact Fuchsian groups $H$ and $K$
accumulates everywhere on the boundary $S^1$; it follows that the map
$\s_\r$ extends by continuity to define a homeomorphism of $S^1$. That
homeomorphism is precisely $\partial\Si_\r$. 

It is now clear that the representation $\Si$ embeds $\Vaut$ in $\Hqs$
as exactly the {\it virtual normalizer of $\Ga$ amongst quasisymmetric
homeomorphisms}. By this we mean that the image by $\Si$ of $\Vaut$ 
is described as:
$$
{\rm Vnorm}_{q.s.}(\Ga) = \{f \in \Hqs: f~ \hbox{conjugates some
finite index}
\leqno(2.13)
$$
$$
~~~~~~~~~~~~~~~~~~~~~~\hbox{subgroup of} ~\Ga~ 
\hbox{to another such subgroup of} ~\Ga \} $$

\noindent
{\bf II.7. $\MCin$ as a subgroup of the universal modular group:}
The representation of $\Vaut$ above allows us to consider 
the action of $\MCin$ on $\Tin$ via the usual type of right 
translations by quasisymmetric homeomorphisms, as is standard 
for the classical action of the universal modular 
group on the universal Teichm\"uller space. 

Recall that the universal Teichm\"uller space of Ahlfors-Bers is the
homogeneous space of right cosets (i.e., $\Mob$ acts by post composition):
$$
\T(1):= {\Hqs}/{\Mob}  
\leqno(2.14)
$$
The coset of $\phi \in \Hqs$ is denoted by $[\phi]$.

Naturally, $\Hqs$ acts as biholomorphic automorphisms of this complex
Banach manifold, $\T(1)$, by right translation (i.e., by pre-composition 
by $f$). In other words, each $f \in \Hqs$ induces the automorphism:
$$
f_{*}:\T(1) \rightarrow \T(1); ~~~ f_{*}([\phi])=[\phi \circ f]
\leqno(2.15)
$$
and this action on $\T(1)$ is classically called the universal 
modular group action (see [N1]).

But having fixed the Fuchsian group $\Gamma$ as above, we see forthwith
from Proposition 2.6 (in II.3) that a copy of the universal 
commensurability Teichm\"uller space, $\Tin$, embeds in $\T(1)$ as follows:

$$
\Tin~~\cong~~\Tin(\Ga) = \{[\phi] \in \T(1): \phi \in \Hqs 
~\hbox{is compatible} 
\leqno(2.16)
$$
$$
~~~~~~~~~~~~~~~~~~~~~~\hbox{with some finite index subgroup of} ~\Ga \}
$$
where the compatibility condition means that there exists some finite
index subgroup $H \subset \Ga$ such that $\phi H {\phi}^{-1} \subset \Mob$. 

\noindent
{\bf Proposition 2.17:} 
{\it The action of $\MCin$ on $\Tin$ coincides with the action, by 
right translations, of the subgroup of the universal modular group 
corresponding to ${\rm Vnorm}_{q.s.}(\Ga) \subset \Hqs$, restricted to 
$\Tin(\Ga) \subset \T(1)$.}

\noindent
{\bf Proof:} By tracing through all the identifications, one finally
needs to verify that for any $f \in {\rm Vnorm}_{q.s.}(\Ga)$, 
the universal modular transformation $f_{*}$ preserves the 
directed union $\Tin(\Ga)$. It is not difficult to verify from the
definition of the virtual normalizer that $f_{*}$ carries  each finite
dimensional stratum in $\Tin(\Ga)$ to another such stratum, and the
Proposition follows.
$\hfill{\Box}$

\noindent
{\bf II.8. Topological transitivity of $\MCin$ on $\Tin$ and allied issues:}
Does $\MCin$ act with dense orbits in $\Tin$? That is a basic query. This
question is directly seen to be equivalent to the following old conjecture 
which, we understand, is due to L.Ehrenpreis and C.L.Siegel:

\noindent
{\bf Conjecture 2.18:} {\it Given any two compact Riemann surfaces, $X_1$
(of genus $g_1 \geq 2$) and $X_2$ (of genus $g_2 \geq 2$), and given any 
$\epsilon > 0$, can one find finite unbranched coverings $\pi_1$ and 
$\pi_2$ (respectively) of the two surfaces such that 
the corresponding covering Riemann surfaces $\tX_1$ and $\tX_2$ are of the 
same genus and there exists a $(1+\epsilon)$ quasiconformal homeomorphism
between them. (Namely, $\tX_1$ and $\tX_2$ come $\epsilon$-close in the 
Teichm\"uller metric.)}

\noindent
{\it Remark:} Since the uniformization theorem guarantees that the {\it
universal} coverings of $X_1$ and $X_2$ are {\it exactly} conformally
equivalent, the conjecture asks whether we can obtain high {\it finite}
coverings that are {\it approximately} conformally equivalent.

\bigskip
\noindent
{\small {\bf III: THE CHARACTERISTIC TOWER AND $\Min$}}
\medskip

\noindent
{\bf III.1. The cofinal set of characteristic covers:}
The unramified finite covering $\pi:\tX \ra X$ is called 
{\it characteristic} if it corresponds to a {\it characteristic 
subgroup} of the fundamental group ${\pi}_{1}(X)$. 
Namely, ${\pi}_{1}(\tX)$ (as a subgroup of ${\pi}_{1}(X)$) must be 
left invariant by every element of $Aut({\pi}_{1}(X))$; this yields 
therefore (by restriction to the subgroup) a homomorphism: 
$$
L_{\pi}~ :~ Aut({\pi}_{1}(X)) ~\longrightarrow ~ Aut({\pi}_{1}(\tX)) 
\leqno(3.1)
$$
Topologically speaking, every diffeomorphism of $X$ lifts to a 
diffeomorphism of $\tX$, and the homomorphism (3.1) corresponds to 
this lifting process. 

Characteristic subgroups are necessarily normal subgroups. 
It is well-known that the normal subgroups of finite index form a 
cofinal family among all subgroups of finite index in $\Ga = 
{\pi}_1(X)$.
We now show the critically important fact that the property continues 
to hold for characteristic subgroups. (Note: All coverings being
considered are {\it finite} and {\it unramified}.)
 
\medskip
\noindent
{\bf Lemma 3.2.}~ {\it The family of finite index characteristic
subgroups, as a directed set partially ordered by inclusion, is
cofinal in the poset of all finite
index subgroups of ${\pi}_{1}(X)$. 
In fact, given any finite covering $f : Y \rightarrow X$, there
exists another finite covering $h : Z \rightarrow Y$ such
that that the composition $f \circ h : Z \rightarrow X$ is a
characteristic cover.}
\medskip

\noindent 
{\bf Proof:} For notational convenience set $G:= {\pi}_1(X)$ 
and $H:={\pi}_1(Y)$, (we will suppress the base points). Using the
monomorphism $\pi_{1}(f)$, the group $H$ will be thought of as a 
subgroup of $G$.

Consider the space of right cosets $S := G/H$, which is a finite
set. The group $G$ has a natural action on $S$ given by the left
multiplication in $G$. So $g\in G$ maps the coset $\{a\}\in S$
to the coset $\{ga\}$. Let $P(S)$ denote the finite group of permutations 
of the set $S$. Let $\rho:G \rightarrow P(S)$ denote
the homomorphism defined by the $G$-action.

Let $\Gamma = Hom(G,P(S))$ denote the set of homomorphisms 
of $G$ into $P(S)$. Since $G$ is a finitely generated group and 
$P(S)$ is a finite group, $\Gamma$ is a finite set.

Define
$$
K =  \bigcap_{\gamma \in \Gamma}{\rm kernel}(\gamma ) \subset G
$$ 
to be the subgroup of $G$ given by the intersection of all the kernels. 
Since $\Gamma$ is a finite set, $K$ is a finite index subgroup of $G$. 
Clearly $K$ is a characteristic subgroup of $G$. 
If we show that $K$ is actually contained in $H$ then the proof 
of the lemma will be complete by taking $h$ to the covering (of $Y$) 
given by the subgroup $K \subset H$.

To prove that $K\subset H$, take any $g \in G$ which is not in $H$,
we will show that $g$ is not in $K$. Consider the action of $g$
on $\{H\}$, the identity coset in $S$. It is mapped to the $\{g\} 
\in P(S)$, the coset given by $g$. Since $g \notin H$, the coset
$\{g\}$ cannot be the coset $\{H\}$. in other words, the
action of $g$ on $P(S)$ is not the trivial action. So $g$ cannot
be in $K$, since $\rho (K) = e$. This completes the proof
of the Lemma.
$\hfill{\Box}$

\smallskip
\noindent
{\bf Alternate proof:} By an argument similar to that used above, 
we see that up to isomorphism there are only finitely many Galois 
coverings of any fixed degree $N$ over a surface $X$ of genus $g$. 

These finitely many normal subgroups of index $N$, sitting within
${\pi}_{1}(X)$, are necessarily permuted amongst themselves by the action
of $Aut({\pi}_{1}(X))$. Taking the intersection of the subgroups
that constitute an orbit under the action of $Aut({\pi}_{1}(X))$
therefore produces a characteristic subgroup of finite index.

As for cofinality, note that any finite index subgroup of any
group $G$ contains within it a subgroup that is normal in $G$ and
is still of finite index. Letting $N$ be the index (in
$G={\pi}_{1}(X)$) of this normal subgroup, and applying the above
construction, we obtain characteristic subgroups of finite index
sitting within any given subgroup of finite index.
$\hfill{\Box}$

\noindent
{\it Example:} Here is a straightforward family of examples for 
finite characteristic coverings of surfaces. Let  
${\pi}_1(X) \rightarrow H_1(X,\Z)$ be the Hurwitz (abelianization) 
map.  Compose this with  the projection  $H_1(X,\Z) \ra H_1(X,{\Z/n})$, 
where $n$ is any integer greater than one. The kernel of this 
composition [${\pi}_1(X)\rightarrow H_1(X,{\Z/n})$] provides
a characteristic subgroup of finite index in $\pi_1(X)$. The 
quotient group, namely the deck transformation group of this 
characteristic covering, is the finite abelian group $H_1(X,{\Z/n})$.

\noindent
{\bf Remarks on fiber-products of coverings}:
Let $f: Y \rightarrow X$ and $g: Z \rightarrow X$ be any 
two pointed coverings of $X$. Let $S$ be the connected 
component of the fiber product $$ S \subset Y\times_X Z$$ containing 
the distinguished point. Let $\mu$ denote the projection of $S$ onto $X$. 
Then the subgroup of $\pi_1(X)$ corresponding to the covering $\mu$ is
simply the {\it intersection} of the two subgroups corresponding to the 
coverings $f$ and $g$. Indeed, if $H$ and $K$ are the subgroups
corresponding to the two given covers, then their fiber product can be
described as the quotient of the universal covering by $H \cap K$.

It follows immediately that {\it any component of the fiber product of two
characteristic coverings over $X$ is also characteristic over $X$.}

Note, of course, that there are factoring projections of $S$ onto 
$Y$ and $Z$ -- denoted by say $\phi$ and $\psi$, respectively. It 
is not in general true that these factoring maps $\phi$ and $\psi$ 
will be characteristic -- even when $f$, $g$ -- and hence $\mu$ -- 
are so. In the definition below of the ordering in the characteristic 
tower over $X$ we are therefore forced to demand that the factoring
morphism should be {\it itself characteristic}. (Otherwise we do not 
get a well-defined inductive system at the moduli spaces level.)

\medskip
\noindent
{\bf III.2. The characteristic tower:}
Consider the tower over the (pointed) surface $X=X_g$ consisting of only 
the {\it characteristic} coverings. Namely, we replace the directed set
(2.4) by the subset: 
$$
\Kch = K^{ch}_g = \{\a \in Mor(\A):\a ~ \hbox{is characteristic and}
~ {\rm Range}(\a)=X \}
$$
For $\a,\b$ in $\Kch$, we say $\b \succ\succ \a$ if and only if 
$\b = \a \circ \theta$ with $\theta$ being also a {\it characteristic}
covering. This gives $\Kch$ the structure of a directed set.

Because of the presence of the homomorphism (3.1), it is evident that
any characteristic cover $\pi$ induces a morphism 
$$
\Mpi: \Mg  \longrightarrow  \Mt
\leqno(3.3)
$$
which is an algebraic morphism between these normal quasi-projective
varieties. In other words,
the map $\Tpi$ of (2.2) {\it descends} to the moduli space level when
the covering $\pi$ is characteristic.

We therefore have a {\it direct system of moduli spaces} over the directed 
set $\Kch$, and passing to the direct limit, we define:
$$
\Min(X) ~ := ~ {\rm ind~lim} \M(X_{g(\a)}), ~~~ \a ~ \in ~ \Kch
\leqno(3.4)
$$
in exact parallel with the definition of $\TinX$ in (2.5). 
(Recall that $X_{g(\a)}$ denotes the domain surface for 
the covering map $\a$.)

\noindent
{\it Question:} Do any two surfaces (genus $g$ and $h$, both
greater than one) have a common characteristic cover? We have been 
unable to resolve this question. Equivalently, we may ask, does 
$\Min(X_{g})$ depend on the genus $g$ of the reference surface? 
Clearly, $\Min(X)$ is naturally isomorphic to $\Min(Y)$ provided 
a common characteristic covering exists.
\medskip

\noindent
{\bf III.3. Mapping-class like elements of $\Vaut$}: If $\a:\tX \ra X$ 
is a morphism of our category $\A$, and $\l:\tX \ra \tX$ is any
self-homeomorphism of $\tX$, then the two-arrow diagram given by 
the two coverings $\a$ and $\a \circ \l$ (the self-correspondence)
defines an element of 
$\Vaut$. Such elements of $\Vaut$ we shall call {\it mapping class like
elements} for obvious reasons (namely, they arise from modular
transformations at some finite covering stage). These elements  
are exactly those virtual automorphisms which {\it fix setwise} some 
finite index subgroup of $\pi_1(X)$. 
We do not know whether every element of $\Vaut$ is mapping class like. 

Utilizing the homomorphisms $L_\a$ of (3.1), we can now define a {\it
direct system of automorphism groups of surfaces} indexed again by
$\Kch$. In fact, we can set:
$$
{\rm Caut}(\pi_{1}(X)) ~ = ~ {\rm dir. lim. Aut}(\pi_{1}(X_{g(\a)})), ~~ 
\a ~ \in ~ \Kch \leqno(3.5)
$$
A little thought shows that the group $\Caut$ consists of those
mapping class like elements which represent automorphisms of finite 
index {\it characteristic} subgroups of $\pi_1(X)$. 

In analogy with the classical situation where $\Mg$ is described as
the quotient of $\Tg$ by the action of the classical mapping class group,
we are now able to describe $\MinX$ in terms of $\TinX$:

\noindent
{\bf Proposition 3.6.}~ {\it  $\Caut$ acts on $\TinX$ to produce
the ind-variety $\MinX$ as the quotient.}

\noindent
{\bf Proof:} Consider the direct system of Teichm\"uller spaces over the
cofinal subset $\Kch$ and let us call $\TchinX$ the corresponding direct
limit space. The inclusion of $\Kch$ in $K(X)$ induces a natural
homeomorphism of $\TchinX$ onto $\TinX$. Clearly, it follows from the 
definition of the group $\Caut$ that $\Caut$ acts on $\TchinX$ to produce
$\MinX$ as the quotient. Therefore, identifying $\TchinX$ with $\TinX$ by
the above homeomorphism, everything follows.  
$\hfill{\Box}$
\medskip

In the paper [BNS] we created determinant bundles over $\TinX$, with the
relating Mumford isomorphisms, the entire construction being invariant
under the full group $\Vaut$. Therefore, in view of the above
Proposition 3.6
it follows immediately that the bundles and isomorphisms constructed in
[BNS] descend to $\MinX$. That is the purport of our main theorem in this
paper, but we shall present the construction independent of the methods in
[BNS]; as we said earlier, our tool in the following chapters will be the 
naturality of the Weil-Petersson K\"ahler forms on the moduli spaces with
respect to the covering maps.

\noindent
{\it A question:} Study the subgroup $\Caut$ in $\Vaut$. Is it a
normal subgroup? Is the index infinite? 

\noindent
{\bf III.4. $\Vaut$ and the Cantor group $\widehat{\pi_{1}(X)}$:}
Consider the algebro-geometric fundamental group of $X$
defined as the profinite completion of the topological fundamental 
group $\Ga=\pi_{1}(X)$. Namely, 
$$
\C\pi_{1}(X) = \widehat{\pi_{1}(X)} = 
{\rm inv~lim} \{\hbox{finite quotients of} ~ \pi_{1}(X)\}
\leqno(3.7)
$$
limit being taken over all finite index normal subgroups of $\Ga$. 
This is the inverse limit of the deck transformation groups of all 
normal (Galois) finite coverings  of $X$. In fact, if we consider the 
inverse limit solenoid construction $\HinX$ running through the cofinal
family of all finite normal covers over $X$, we see that the fiber of the
fibration $\pi_{\infty}:H_{\infty}(X) \rightarrow X$ 
is precisely this Cantor-set group $ \widehat{\pi_{1}(X)}$. It is not 
hard to see that there is a natural embedding of $\Vaut$ into the virtual 
automorphism group of this Cantor group. Regarding this relationship, 
and concomitant matters, we will have more to say in a forthcoming 
article [NaSa].

\bigskip
\noindent
{\small {\bf IV: CURVATURE FORMS OF $DET$ BUNDLES ON $\Mg$}}
\medskip

\noindent 
{\bf IV.1. Line bundles on the moduli space:}
There are several closely related concepts of line bundles associated
to the moduli spaces of Riemann surfaces. We will recall the definition
of the Picard group $\Pf$ -- which is the most basic one from the 
algebro-geometric standpoint. 
$\Pf$ denotes {\it the Picard group of the moduli functor}. 
An element of $\Pf$ consists in prescribing an algebraic line bundle
$L_F$ on the base space $S$ for every algebraic family 
$F=(\gamma: V \rightarrow S)$ of Riemann surfaces of genus $g$ over 
any quasi-projective base $S$. Moreover, for every commutative diagram of
families $F_{1}$ and $F_{2}$ having the morphism $\a$ from the base
$S_{1}$ to $S_{2}$, there must be assigned a corresponding isomorphism
between the line bundle $L_{F_{1}}$ and the pullback via $\a$ of the bundle
$L_{F_{2}}$. For compositions of such pullbacks, these isomorphisms between
the prescribed bundles must satisfy the self-evident compatibility 
condition. Two such prescriptions of line bundles over bases $S$
define the same element of $\Pf$ if there are compatible isomorphisms
between the bundles assigned for each $S$. See [Mum], [HM], [AC] for
details. [Note: Mumford has considered this Picard group of the moduli
functor also over the Deligne-Mumford compactification of $\Mg$.]

\noindent
{\it The Hodge line bundle:}
We introduce this fundamental (generating!) element of $\Pf$. Consider 
any smooth family of genus $g$ Riemann surfaces, 
$F:=(\gamma:E\rightarrow S)$. The ``Hodge bundle'' on the 
parameter space $S$ is defined to be the dual of $\ext{g}(R^1{\ga}_*{\O})$. 
Here the $R^1$ denote the usual first direct image (see, for example, [H]). 
Associating to each family $F$ its Hodge line bundle, one 
obtains an element of $\Pf$, per definition. 
The fiber of the Hodge line bundle over the point $s \in S$ is the
top exterior product $\ext{g}H^1({X_s},{\O})^{*}$, where $X_s$ denotes the
genus $g$ curve ${\ga}^{-1}(s)$. By the Serre duality, this exterior
product is canonically isomorphic to $\ext{g}H^0(X_s, K)$, where
$K=K_{X_s}$ is its cotangent bundle. It is  a fundamental fact
that $\Pf$ is generated by the Hodge line bundle [AC]. Moreover, for
$g\geq 3$, the group $\Pf$ is freely generated by the Hodge bundle.
In particular, for $g\geq 3$, we have $\Pf = \Z$. (For $g=2$, 
$\Pf = \Z/\!{10\Z}$.)

\smallskip
\noindent 
{\bf The relation between ${\rm Pic}(\Mg)$ and ${\rm Pic_{hol}}(\Mg)$
with $\Pf$:}
Let $DET_0 \longrightarrow {\cal T}_g$
be the Hodge bundle on the Teichm\"uller space. There is a natural
lift of the action of the modular group $MC_g$ on ${\cal T}_g$ to
$DET_0$. Assume that $g\geq 2$. Since the automorphism group of a Riemann
surface of genus at least two is a finite group, there is a positive 
integer $n(g)$, (for example, $[84(g-1)]!$ works) such that the 
induced action of any isotropy subgroup for the action of $MC_g$ 
on ${\cal T}_g$, on the fiber of $DET^{m.n(g)}_0$, for any $m \in \Z$, 
is the trivial action. Consequently, each of the line bundles
$DET^{m.n(g)}_0$ descends as an algebraic line bundle on $\Mg$. All 
algebraic line bundles on $\Mg$ are known to arise this way. 

The Picard group of $\Mg$, denoted by ${\rm Pic}(\Mg)$, consisting 
of isomorphism classes of algebraic line bundles on $\Mg$, is a finite 
index subgroup of $\Pf$ -- see [AC]. Any holomorphic line bundle on 
the Teichm\"uller space ${\cal T}_g$, equipped with a lift of the 
action of the mapping class group $MC_g$, such that the
action of the isotropy subgroup of any point on the fiber is trivial,
must be a power of the Hodge line bundle for the universal family of 
Riemann surfaces over ${\cal T}_g$. Let ${\rm Pic_{hol}}(\Mg)$ denote 
the group of isomorphism classes of holomorphic line bundles $\Mg$. Then 
from the above remarks it follows that we have (for $g \geq 3$):
$$
{\rm Pic}(\Mg)\otimes_{\Z}\Q ~ = ~ {\rm Pic_{hol}}(\Mg)\otimes_{\Z}\Q
~ = \Pf\otimes_{\Z}\Q ~ = ~ \Q
\leqno(4.1)
$$

\smallskip
\noindent
{\bf $DET$ bundles for families:} Given, as before, 
any Kodaira-Spencer family $F=(\gamma:V \rightarrow S)$,
of compact Riemann surfaces of genus $g$, and a holomorphic vector 
bundle $E$ over the total space $V$, we can consider the base $S$ as
parametrizing a family of elliptic d-bar operators.
The operator corresponding to $s \in S$ acts along the fiber Riemann
surface $X_s = {\gamma}^{-1}(s)$ :
$$
\bar{\partial}_{s}: C^{\infty}({\gamma}^{-1}(s), E) ~ \longrightarrow ~
C^{\infty}({\gamma}^{-1}(s), E \otimes {{\Omega}^{0,1}_{X_{s}}})
$$
One defines the associated vector space of one dimension given by:
$$
DET(\bar{\partial}_{s}) = (\ext{\rm top}{\rm ker}\bar{\partial}_{s}) 
\otimes (\ext{\rm top}{\rm coker}\bar{\partial}_{s})^{*}
\leqno(4.2)
$$
and it is known that these complex lines fit together naturally over the
base space $S$ giving rise to a holomorphic line bundle over $S$ called 
$DET(\bar{\partial})$. In fact, this entire 
``determinant of cohomology'' construction is natural with
respect to morphisms of families and pullbacks of vector bundles. Note
that the definition of the determinant line in (4.2) coincides with that
given in [D], but is dual to the one in [Bos].

We could have followed the above construction through for 
the {\it universal genus $g$ family $V_g$ over $\Tg$} 
(see [N1]), with the vector bundle $E$ being, variously, 
the trivial line bundle over the universal curve, or the 
vertical (relative) tangent bundle, or any of its tensor powers.
It is easy to verify the following: setting $E$ to be the 
trivial line bundle over $V$ for any family 
$F=(\gamma:V \rightarrow S)$, the above prescription
for $DET$ provides merely another description of the Hodge line bundle.

By the same token, setting over any family $F$ the vector bundle $E$ to be 
the $m^{th}$ tensor power of the vertical cotangent bundle along the 
fibers, we get by the $DET$ construction a well-defined member 
$$
DET_{m}={\l}_{m} ~ \in ~ \Pf, ~~~~ m \in \Z.
\leqno(4.3)
$$
Serre duality shows that $DET_{m}=DET_{1-m}$, in $\Pf$.
Clearly, $\l_{0}$ is the Hodge bundle, and by ``Teichm\"uller's 
lemma'' (see [N1]) one notes that ${\l}_{2}$ represents 
the canonical bundle of the moduli space; indeed, the fiber
of $DET_{2}$ at any Riemann surface $X \in \Mg$ is the top 
exterior product of the space of holomorphic quadratic 
differentials on $X$.

\medskip
\noindent
{\bf IV.2. Mumford isomorphisms:} By applying the Grothendieck-Riemann-Roch
theorem it was proved by Mumford in [Mum] that as elements of $\Pf$ one has 
$$
{\l}_m ~ = ~ (6m^{2} - 6m + 1)\hbox{-th} ~~ \hbox{tensor power of Hodge}
~(=\l_{0}) 
\leqno(4.4)
$$

The complement of $\Mg$ in its Satake compactification is of 
codimension at least two if $g\geq 3$. The Hartogs theorem implies
that there are no non-constant holomorphic functions on $\Mg$ ($g \geq 3$).
Therefore the choice of an isomorphism of ${\l}_{m}$ with 
${{\l}_{0}}^{\otimes {(6m^{2}-6m+1)}}$ is {\it unique} up to a nonzero
scalar. We would like to put canonical hermitian metrics on these DET
bundles so that this essentially unique isomorphism actually becomes an
unitary isometry. This follows from the theory of the:

\medskip
\noindent 
{\bf IV.3. Quillen metrics on DET bundles:} 
If we prescribe a conformal Riemannian metric on the fiber Riemann
surface $X_s$, and simultaneously a hermitian fiber metric on the vector
bundle $E_s$, then clearly this will induce a natural $L^2$ pairing on the
one dimensional space $DET(\bar{\partial}_{s})$ described in (4.2). Even if
one takes a smoothly varying family of conformal Riemannian metrics on the
fibers of the family, and a smooth hermitian metric on the vector bundle
$E$ over $V$, these $L^2$ norms on the DET-lines may fail to fit together
smoothly (basically because the dimensions of the kernel or cokernel for 
${{\bar{\partial}}_{s}}$ can jump as $s$ varies over $S$). However, 
Quillen, and later Bismut-Freed and other authors, have described
a ``Quillen modification'' of the $L^2$ pairing which always produces a 
smooth Hermitian metric on $DET$ over $S$, and has important functorial 
properties.  

\noindent
{\bf Remark:}
Actually, in the cases of our interest the usual Riemann-Roch 
theorem shows that the dimensions of the kernel and cokernel spaces 
remain constant as we vary over moduli -- so that the $L^2$ metric 
is itself smooth. Nevertheless, the Quillen metric will be crucially 
utilized by us because of certain functorial properties, and curvature 
properties, that it enjoys.
\medskip

Using the metrics assigned on the Riemann surfaces (the fibers of
$\gamma$), and the metric on $E$, one gets $L^2$ structure on the spaces of
$C^{\infty}$ sections that constitute the domain and target for our d-bar
operators. Hence ${\bar{\partial}}_{s}$ is provided with an adjoint
operator ${{\bar{\partial}}_{s}}^{*}$, and one can therefore construct the
positive (Laplacian) elliptic operator as the composition:
$$
{\Delta}_{s} ~ = ~ {{\bar{\partial}}_{s}}^{*} \circ {{\bar{\partial}}_{s}}, 
$$
These Laplacians have a well-defined (zeta-function regularized)
determinant, and one sets:
$$
\hbox{Quillen norm on fiber of}~ DET ~=~ 
(L^{2}~\hbox{norm on that fiber}).({\rm det}{\Delta}_{s})^{-1/2}
\leqno(4.5)
$$
{\it This turns out to be a smooth metric on the line bundle $DET$.} 
See [D], [Q], [BF], [BGS].

In the situation of our interest, the vector bundle $E$ is the vertical
tangent (or cotangent) line bundle along the fibers of $\gamma$, or its
powers, so that the assignment of a metric on the Riemann surfaces already
suffices to induce a Hermitian metric on $E$. Hence one gets a Quillen norm
on the various $DET$ bundles ${\l}_m$ ($\in \Pf$) for every choice
of a smooth family of conformal metrics on the Riemann surfaces.
{\it The Mumford isomorphisms (over any base $S$) become isometric
isomorphisms with respect to the Quillen metrics.}

Let $T_{\rm vert} \rightarrow V$ denote the relative tangent bundle, namely
the kernel of the differential map of the projection of $V$ onto $S$. 
The {\it curvature form} (i.e., first Chern form) on the base 
$S$ of the Quillen DET bundles has a particularly elegant expression: 
$$
c_{1}(DET, Quillen~metric) ~ = ~ \int_{V|S}^{}(Ch(E) Todd({T}_{\rm vert})) 
\leqno(4.6)
$$
where the integration represents integration of differential forms along
the fibers of the family $\gamma: V \rightarrow S$ [Bos], [D].

We now come to one of our main tools in this paper.
By utilizing the uniformization theorem (with moduli parameters), the
universal family of Riemann surfaces over $\Tg$, and hence any holomorphic
family $F$ as above, has a natural smoothly varying family of Riemannian 
metrics on the fibers given by the constant curvature $-1$ {\it Poincar\'e
metrics}. The Quillen metrics arising on the $DET$ bundles $\l_m$ from
the Poincar\'e metrics on $X_s$ has the following fundamental property 
for its curvature:
$$
c_1({\l}_{m}, Quillen) ~ = ~{\frac{1}{12\pi^{2}}}(6m^{2}-6m+1)
{{\omega}_{WP}},
\hspace{.5in} m \in \Z
\leqno(4.7)
$$
\noindent
where ${\omega}_{WP}$ denotes the (1,1) K\"ahler form on $\Tg$ for the
classical {\it Weil-Petersson} metric of $\Tg$. We remind the reader that
the cotangent space to the Teichm\"uller space at $X$ can be canonically
identified with the vector space of holomorphic quadratic differentials
on $X$, and the WP Hermitian pairing is obtained as 
$$
(\phi, \psi)_{WP} = \int_{X}^{} \phi {\bar {\psi}} (Poin)^{-1}
\leqno(4.8)
$$ 
Here $(Poin)$ denotes the area form on $X$ induced by the Poincar\'e
metric. That the curvature formula (4.6) takes the special form (4.7) for
the Poincar\'e family of metrics has been shown by Wolpert [Wol] and 
Zograf-Takhtadzhyan [ZT].

Indeed, (4.6) specialized to
$E=T_{\rm vert}^{\otimes{-m}}$ becomes simply $(6m^{2}-6m+1)/12$ times
${\int_{V|S}^{}}{c_{1}{(T_{\rm vert})}^{2}}$. This last integral
represents, for the Poincar\'e-metrics family, ${\pi}^{-2}$ times
the Weil-Petersson symplectic form. See also [BF], [BGS], [BK], 
[Bos], [Wol], [ZT].

Applying the above machinery, we will investigate the behaviour of 
the Mumford isomorphisms in the situation of a covering map between 
surfaces of different genera.

\bigskip
\noindent 
{\small{\bf V. CHARACTERISTIC COVERINGS AND DET BUNDLES:}}
\medskip

\noindent
{\bf V.1. Comparison of Hodge bundles:} 
Let $\pi: \tX \rightarrow X$ be a {\it characteristic 
covering} of degree $N$. Recall from topology that
$N=(\tg-1)/(g-1)$, where $\tg$ and $g(\geq 2)$ are respectively the
genera of $\tX$ and $X$.  Let ${\cal M}(\pi):\Mg  \rightarrow  \Mt$
be the morphism induced by $\pi$ as in (3.3). We are now 
in a position to compare the two candidate Hodge bundles
that we get over $\Mg$ -- one is the pullback of the Hodge 
bundle from $\Mt$ using ${\cal M}(\pi)$, and the other being 
the Hodge bundle of $\Mg$ itself. The same comparison
will be worked out simultaneously for all the $DET_m$ bundles. 

\noindent 
{\it Notations:}
Let $\l = {DET}_{0}$ denote, as before, the Hodge bundle on $\Mg$
(a member of $\Pf$, as explained), and let ${\tilde{\lambda}}$ 
denote the Hodge line bundle over $\Mt$. Further, let $\om = {\om}_{WP}$ 
and $\tom$ represent the Weil-Petersson forms (i.e., the K\"ahler 
forms corresponding to the WP Hermitian metrics) on $\Mg$ and $\Mt$, 
respectively. 

The naturality of Weil-Petersson forms under coverings is manifest in the
following basic Lemma:

\medskip
\noindent 
{\bf Lemma 5.1.} {\it The 2-forms $N(\om)$ and $({\Mpi})^{*}\tom$ 
on $\Mg$ coincide.}
\medskip

\noindent {\bf Proof:}~ This is basically a straightforward
computation. Recall that the cotangent space to the Teichm\"uller 
space is canonically isomorphic to the space of quadratic 
differentials for the Riemann surface represented by that point.
(It is actually sufficient to prove this Lemma at the Teichm\"uller level.)
Now, at any point $\a \in \Mg$, the co-derivative morphism is a map 
on cotangent spaces induced by the map $\Mpi$:
$$
(d{\Mpi})^*: T^*_{\Mpi(\a)}\Mt \longrightarrow T^*_{\a}\Mg
$$
The action of this map on any quadratic differential $\p$ on the 
Riemann surface $\Mpi(\a)$, (i.e., $\p \in H^0(\Mpi(\a), K^2)$),  
is given by:
$$
(d{\Mpi})^*{\p}= 1/N(\sum_{f\in Deck}^{}f^*{\p})
\leqno(5.1)
$$ 
Here $Deck$ denotes, of course, the group of deck transformations
for the covering $\pi$. Now recall that a covering map $\pi$ induces a
local isometry between the respective Poincar\'e metrics, and that $N$ 
copies of $X$ will fit together to constitute $\tX$. The lemma 
therefore follows by applying formula (5.1) to two quadratic 
differentials, and pairing them by the Weil-Petersson pairing as 
per definition (4.8). 
$\hfill{\Box}$

\smallskip
The above Lemma, combined with the fact that the curvature form of the 
Hodge bundle is $(12{\pi}^{2})^{-1}$ times $\om_{WP}$, shows that the
{\it curvature forms (with respect to the Quillen metrics) of the 
two line bundles ${\l}^{\otimes N}$ and $({\Mpi})^*{\tilde{\lambda}}$ 
coincide}. Do the bundles themselves coincide? Yes: 

\smallskip
\noindent 
{\bf Theorem 5.2:} {\it  
Let $\pi: \tX \longrightarrow X$ be a characteristic covering of 
degree $N$. Then the two line-bundles ${\l}^{\otimes N}$ and 
$({\Mpi})^*{\tilde{\lambda}}$, as members of $\Pf$, are equal.

If $g\geq 3$, then such an isomorphism:
$$
F_{\pi} :{\l}^{\otimes N}  \longrightarrow  ({\Mpi})^*{\lt}
$$
is uniquely specified up to the choice of a nonzero scaling constant.
Up to a constant any such isomorphism must be an unitary isometry
between the $d$-th power of the Quillen metric on $\l$ and the Quillen
metric of $\lt$.

The same assertions hold for each of the bundles $DET_m$, 
$m=0,1,2,\cdots$. Namely, the bundles 
${DET_m}^{\otimes N}$  and $({\Mpi})^*{\widetilde{DET_m}}$
are isometrically isomorphic. (By ${\widetilde{DET_m}}$ we mean,
of course, the $DET_m$ construction over $\Mt$.)}

\medskip

\noindent {\bf Proof:} The basic principle is that ``curvature of the 
bundle determines the bundle uniquely'' over  any space whose Picard
group is discrete, in particular, therefore over 
the moduli space $\Mg$. That is automatic since the first 
Chern class is given by curvature form, and it is known that 
$c_1: \Pa \rightarrow H^{2}(\Mg,\Z)=\Z$ is injective.

The uniqueness assertion follows since two isomorphisms can only differ by
a global holomorphic function on the base space. But, as explained in 
Section IV.2, $\Mg$ does not admit any 
nonconstant global holomorphic functions for $g \geq 3$.

Consider now the pull-back of the Quillen metric on $\lt$ to ${\l}^{N}$. 
We know that the curvature form coincides with the 
curvature of the $N$-th power of the Quillen
metric of $\l$. Therefore the ratio of these two metrics on ${\l}^{N}$ is a
positive function $f$ on the base, with $\log(f)$ pluri-harmonic. But 
there is no obstruction in making $\log(f)$ the real part of a global 
holomorphic function on $\Mg$ since the first Betti number of 
$\Mg$ is zero. By the Satake-Hartogs argument then, $f$ must be constant, 
implying the isometrical nature (up to constant) of any isomorphism, as
desired.
$\hfill{\Box}$

\smallskip
A characteristic covering $\pi$ of degree $N$ therefore 
provides a map from ${\rm Pic_{fun}}(\Mt)$ into $\Pf$ that is an
embedding of infinite cyclic groups with index equal to the 
degree $N$ of the covering. Indeed, the above result proves that the
Hodge bundle on the smaller moduli space, raised to the tensor power $N$,
extends over the larger moduli space as the Hodge bundle thereon. 

Since we will deal with towers of characteristic coverings, we need 
to make sure that the isomorphisms we pick up are perfectly 
{\it compatible} -- in order to create, as appropriate, inductive or
projective systems of objects by utilizing these isomorphisms. 

Therefore let $\nu:Z \rightarrow X$ be a characteristic covering of
degree $N$ of a surface $X$ of genus $g (\geq 3)$, and suppose that $\nu$ 
allows a decomposition into a pair of {\it characteristic} coverings, 
$f_1$ and $f_2$ of orders $n_1$ and $n_2$ respectively. Namely,
$f_2 \circ  f_1 = \nu$ and $N={n_1}{n_2}$,
$$
\matrix{Z&\mapright{f_1}&Y&\mapright{f_2} &X\cr}
$$
Denote the genera of $Y$ and $Z$ by $g_1$ and $g_2$ respectively. 
Let 
$\M(f_{1}): {\cal M}_{g_1}\rightarrow\,{\cal M}_{g_2}$,
$\M(f_{2}):\Mg \rightarrow\,{\cal M}_{g_1}$ and
$\M(\nu) = \M(f_2\circ f_1):\Mg \rightarrow\,{\cal M}_{g_2}$
be the induced maps of moduli spaces as in (3.3).
We let ${\l}_1$ and ${\l}_2$ denote the Hodge bundles on
${\cal M}_{g_1}$ and ${\cal M}_{g_2}$, respectively, and 
$\l$ the Hodge on $\Mg$. The compatibility we need to establish 
is that the following diagram {\it commutes}:
 
\smallskip
$$
\matrix{{\l}^{n_1n_2}&\mapright{id}&{\l}^{n_1n_2}\cr
\mapdown{F_{f_2}^{\otimes{n_1}}}&&\mapdown{F_\nu}\cr
{\M(f_{2})}^*{{\l}_{1}^{n_1}}&\mapright{{\M(f_{2})}^{*}(F_{f_1})}&
{\M(\nu)}^*{\l}_{2}\cr}
$$
In the diagram above, the morphisms $F_\nu$, $F_{f_1}$ and $F_{f_2}$ 
are obtained as instances of the map $F_\pi$ of Theorem 5.2, applied 
when $\pi$ is taken in turn to be each of the three coverings under 
scrutiny. 
 
But Theorem 5.2 asserts that both the mappings
$F_{\nu}$ and $(\M({f_2}))^{*}(F_{f_1}) \circ {F_{f_2}}^{\otimes{n_1}}$ 
represent isomorphisms between the bundles ${\l}^{N}$ and 
${\M(\nu)}^*{\l}_{2}$, over $\Mg$ -- hence they can only differ by a
global holomorphic function on the base space $\Mg$. By the Satake-Hartogs 
argument, recall that any holomorphic function on $\Mg$ must be a 
constant, so we have proved that the two mappings above are necessarily 
just scalar multiples of each other.
Moreover, since the morphism $F_\pi$ in Theorem 5.2 was chosen to be 
an isometry with respect to the appropriate powers of Quillen metrics, 
the scalar under concern must be of norm one. By adjusting $F_\nu$
by the appropriate scale factor (which is anyway at our disposal), 
we can therefore choose the morphisms involved so that the diagram 
commutes, as desired.

\medskip
\noindent
{\bf V.2. The Main Theorems:}
We are now in a position to formulate our main results as the
construction of a certain sequence of canonical $DET_m$ line bundles 
over the ind-spaces ${\TchinX}$ and $\Min$ and obtain the desired 
Mumford isomorphisms between the relevant tensor powers of these bundles. 

\noindent
{\bf Line bundles on ``ind-spaces'':} 
A line bundle on the inductive limit of an inductive system of varieties or
analytic spaces, is, by definition ([Sha]), a collection of line bundles 
on each stratum (i.e., each member of the inductive system of spaces) 
together with compatible bundle maps. The compatibility condition for 
the bundle maps is the obvious one relating to their behavior with 
respect to compositions, and guarantees that the
bundles themselves fit into an inductive system.

Now, a line bundle with Hermitian metric on an inductive limit space 
is a collection of hermitian metrics for the line bundles over each 
stratum such that the connecting bundle maps are unitary. 
The isomorphism class of such a direct system of Hermitian line bundles 
(over a direct system of spaces), can clearly be thought
of as an element of the inverse limit of the groups consisting of
isomorphism classes of holomorphic Hermitian line bundles on the
stratifying spaces. (The group operation is defined by tensor product.)

For any complex space $M$, let us denote by ${\rm Pich}(M)$ the group
consisting of the isometric isomorphism classes of holomorphic Hermitian 
line bundles on $M$. Moreover, let ${\rm Pich}(M)_{\Q}$ denote 
${\rm Pich}(M)\otimes_{\Z}\Q$, this is constituted 
by the isomorphism classes of ``rational'' holomorphic 
Hermitian line bundles over $M$.  

For any {\it inductive} system of spaces $M_i$, one obtains a 
corresponding {\it projective} system of groups ${\rm Pich}(M_{i})$ -- 
whose limit will be denoted by 
$\lim_{\leftarrow}{\rm Pich}(M_{i})$.
A rational Hermitian line bundle over the inductive limit space 
$\lim_{\rightarrow}{M_{i}}$ is then, by definition, an element of 
$\lim_{\leftarrow}{\rm Pich}(M_{i})_{\Q}$.

Our main result is to create natural elements, related by the relevant
Mumford isomorphisms, of 
$\lim_{\leftarrow}{\rm Pich}({\cal T}_{g_i})_{\Q}$ (and of 
$\lim_{\leftarrow}{\rm Pich}({\cal M}_{g_i})_{\Q}$), as we go through 
the directed tower of all characteristic coverings over a fixed base 
surface $X$ of genus $g$.  

\smallskip
\noindent 
{\bf Theorem 5.3: ``Universal $DET$ line bundles'':}~
{\it There exist canonical elements of the inverse limit
$\lim_{\leftarrow}{\rm Pich}({\cal T}_{g_i})_{\Q}$, namely 
hermitian line bundles on the ind-space $\TchinX$, 
representing the Hodge and higher DET bundles with
respective Quillen metrics:    
$$
{\La}_{m}~ \in ~ \lim_{\leftarrow}{\rm Pich}({\cal T}_{g_i})_{\Q},
\hspace{.5in} m \in \Z
$$
The pullback (namely restriction) of $\La_m$ to each of the stratifying
Teichm\"uller spaces ${\cal T}_{g_i}$ is $(n_i)^{-1}$ times the 
corresponding determinant bundle $DET_{m}$, (with $(n_i)^{-1}$ times 
its Quillen metric), over the stratum ${\cal T}_{g_i}$. Here $n_i$
denotes the degree of the covering from genus $g_i$ to genus $g$.}

{\it Precisely the same statements as above go through on the
ind-variety $\MinX$. The universal line bundles ${\La}_m$ live on 
$\MinX$, and restrict to each stratum ${\cal M}_{g_i}$ as 
$(n_i)^{-1}$ times the corresponding ${\l}_{m}$ bundle 
over ${\cal M}_{g_i}$.}

\noindent
{\bf Proof:} We may work with modular-invariant bundles over the
Teichm\"uller spaces, the construction over the inductive limit of
moduli spaces being identical. 

The foundational work is already done in Section V.1 above. In fact,
let ${\l}_{0,i}$ 
represent the Hodge bundle with Quillen
metric in ${\rm Pich}({\cal T}_{g_i})_{\Q}$. Then, for any $i\in I$, 
taking the element 
$$
(1/{n_i}){\l}_{0,i} ~ \in ~ {\rm Pich}({\cal T}_{g_i})_{\Q}
$$ 
provides us a compatible family of hermitian line bundles (in the rational
Pic) over the stratifying Teichm\"uller spaces -- as required in the
definition of line bundles over ind-spaces. The connecting family of 
bundle maps is determined (up to a scalar) by Theorem 5.2. 

Notice that prescribing a base point in $\Tg$ fixes a compatible
family of base points in each Teichm\"uller space ${\T}_{g_i}$ 
(and, therefore, also in each moduli space ${\cal M}_{g_i}$). 
If we  choose a vector of unit norm in the fiber 
over each of these base points, then that procedure 
rigidifies uniquely all the scaling factor ambiguities 
in the choice of the connecting bundle maps. Then the connecting 
unitary bundle maps for the above collection become {\it compatible}, 
and we have therefore constructed the universal Hodge, $\La_{0}$, over 
${\TchinX}$ (and, by the same proof, over $\MinX$).

Naturally, the above analysis can be repeated verbatim for
each of the d-bar families, and one thus obtains elements 
$\La_{m}$ for each integer $m$. Again the pullback of 
$\La_{m}$ to any of the stratifying ${\cal T}_{g_i}$ produces
$(n_i)^{-1}$ times the $m$-th $DET$ bundle (with appropriate power
of the Quillen metric) living over that space. 
$\hfill{\Box}$

\smallskip
\noindent 
{\bf Theorem 5.4: ``Universal Mumford isomorphisms'':}~
{\it Over the inductive limit of Teichm\"uller space ${\TchinX}$,
or over the ind-variety $\MinX$, we have for each $m \in \Z$:
$$
\La_{m} ~ = ~ {\La_{0}}^{\otimes (6m^{2}-6m+1)}
$$ 
as equality of hermitian line bundles.} 
\smallskip

\noindent 
{\bf Proof:} Follows directly from the genus-by-genus isomorphisms
of (4.4) and our universal line bundle construction above.
$\hfill{\Box}$

\newpage

 
\medskip
\noindent
{\underline{Authors' Addresses:}}

\medskip
\noindent
Tata Institute of Fundamental Research,\\
Colaba, Bombay 400 005, India; ``indranil@math.tifr.res.in''\\
{\it and:} Institut Fourier, 38402 Saint Martin d'Heres, France
 
\medskip
\noindent
The Institute of Mathematical Sciences,\\
CIT Campus, Madras 600 113, India; ``nag@imsc.ernet.in''\\
{\it and:} University of Southern California, Los Angeles,\\ 
California 90089-1113, USA; ``nag@math.usc.edu'' 

\end{document}